\documentclass[10pt]{article}
\usepackage{amssymb, amsmath, url, graphicx, setspace, geometry}
\usepackage{calrsfs}
\usepackage{wasysym, xcolor}

\DeclareMathOperator{\os}{os}

\title{\bf On two problems about order sequences of finite groups}
\author{Mihai-Silviu Lazorec}
\date{October 27, 2024}

\begin{document}

\maketitle

\begin{abstract}
The order sequence of a finite group $G$ is a non-decreasing finite sequence formed of the element orders of $G$. Several properties of order sequences were studied by P. J. Cameron and H. K. Dey in a recent paper that concludes with a list of open problems. In this paper we solve two of these problems by showing the following facts: 1) if there is a non-supersolvable/non-solvable group of order $n$, it is not always true that its order sequence is properly dominated by the order sequence of any supersolvable/solvable group of order $n$; 2) the supersolvability of a finite group cannot be described by its order sequence.  
\end{abstract}



\section{Introduction}

Let $G$ be a finite group of order $n$. We denote the order of an element $g$ of $G$ by $o(g)$. By placing the element orders of $G$ in non-decreasing order, we get the so-called order sequence of $G$. We denote this sequence by $\os(G)$. Formally, if $G=\{ g_1, g_2, \ldots, g_n\}$, then the order sequence of $G$ is
$$\os(G)=(o(g_1), o(g_2), \ldots, o(g_n)), \text{ \ where \ } o(g_i)\leq o(g_{i+1}), \ \forall \ i \in \{ 1, 2, \ldots, n-1\}.$$
For the ease of writing, when the order sequence includes several repeated elements, we write some pairs such as $(o(g), i)$ into $\os(G)$ meaning that the element order $o(g)$ has a multiplicity of $i$ (i.e. $o(g)$ appears in $\os(G)$ for $i$ times). For example, for the cyclic group $C_n$ of order $n$, we could write
$$\os(C_n)=((1, 1), (d_1, \varphi(d_1)), (d_2, \varphi(d_2)), \ldots, (d_j, \varphi(d_j)), (n, \varphi(n))),$$ where $d_1, d_2, \ldots, d_j$ are the non-trivial divisors of $n$ and $\varphi$ is Euler's totient function. 

The origin of this paper is related to \cite{5}, where P. J. Cameron and H. K. Dey introduced the order sequence of a finite group and studied several problems concerning it. In what follows, we briefly recall some ideas and results from their paper. If $G=\{ g_1, g_2, \ldots, g_n\}$ and $H=\{ h_1, h_2, \ldots, h_n\}$ are two groups of order $n$, we say that $\os(G)$ dominates $\os(H)$ if $o(g_i)\geq o(h_i)$, for all $i\in\{ 1, 2, \ldots, n\}$. This domination relation is a partial order on the set of order sequences of groups of order $n$. The order sequences of the alternating group $A_4$ and of the dihedral group $D_{12}$ are
$$\os(A_4)=(1, 2, 2, 2, 3, 3, 3, 3, 3, 3, 3, 3) \text{ \ and \ } \os(D_{12})=(1, 2, 2, 2, 2, 2, 2, 2, 3, 3, 6, 6),$$
showing that the domination relation is not a total order. If $K=\{ k_1, k_2, \ldots, k_m\}$ is another group of order $m$, then $\os(G)\os(K)$ is a sequence formed of all products $o(g_i)o(k_j)$ arranged in non-decreasing order. The sequence $\os(G)\os(K)$ is not necessarily an order sequence of a finite group. Regarding this aspect, the authors proved the following results.\\

\textbf{Theorem 1.1.} \textit{Let $G_1, G_2$ and $H$ be finite groups. Then
\begin{itemize}
\item[i)] $\os(G_1\times G_2)=\os(G_1)\os(G_2)$ if and only if $G_1$ and $G_2$ are of coprime orders;
\item[ii)] If $G_1$ and $G_2$ are not of coprime orders, there is no finite group $K$ such that $\os(K)=\os(G_1)\os(G_2);$
\item[iii)] If $|G_1|=|G_2|$, $\os(G_1)$ dominates $\os(G_2)$ and $G_1$ and $H$ are of coprime orders, then $\os(H\times G_1)$ dominates $\os(H\times G_2)$.
\end{itemize}}     

Another result proved in the same paper concerns the domination relation with respect to finite groups' nature. More exactly, the following theorem holds.\\
  
\textbf{Theorem 1.2.} \textit{Let $n$ be a positive integer such that there exists a non-nilpotent group $H$ of order $n$. Then the order sequence of any nilpotent group $G$ of order $n$ properly dominates $\os(H)$.}\\

The original form of one of the several open questions posed by the authors in the last section of \cite{5} can be stated as follows.\\

\textbf{Question 1.3.} \textit{Let $n$ be a positive integer such that there exists a non-solvable group $H$ of order $n$. Is is true that the order sequence of any solvable group $G$ of order $n$ properly dominates $\os(H)$?}\\

The first main objective of our paper is to provide a negative answer to the above question. We are going to show that there are infinitely many possible choices for $n$ such that there is a unique non-solvable group $H$ of order $n$ which properly dominates the order sequence of at least one solvable group $G$ of order $n$. Moreover, we study the 'supersolvable' version of the same question and, again, we prove that there are infinitely many possible choices for $n$ such that there is a unique non-supersolvable group $H$ of order $n$ which properly dominates the order sequence of at least one supersolvable group $G$ of order $n$. \\ 

The cyclicity of a finite group can be described by its order sequence. More exactly, if $G$ is a finite group and $\os(G)=\os(C_{|G|})$, then $G\cong C_{|G|}$, so $G$ is cyclic. Indeed, since $\os(G)=\os(C_{|G|})$, it follows that $\psi(G)=\psi(C_{|G|})$, where $\psi(G)$ and $\psi(C_{|G|})$ are  the sums of element orders of $G$ and $C_{|G|}$, respectively. Then, the main result of \cite{2} guarantees that $G\cong C_{|G|}$. However, the commutativity of a finite group cannot be described by its order sequence. For instance, for an odd prime $p$, $C_p^3$ has the same order sequence as the non-abelian group
$$He_p=\langle x, y, z \ | \ x^p=y^p=z^p=1, [x, z]=[y, z]=1, [x, y]=z\rangle,$$
whose exponent is $p$. Going back to \cite{5}, the authors proved that the nilpotency of a finite group can be characterized by its order sequence.\\

\textbf{Theorem 1.4.} \textit{Let $G$ and $H$ be finite groups. If $G$ is nilpotent and $\os(G)=\os(H)$, then $H$ is nilpotent.}\\

Another open question posed by the authors concerns the possibility of describing the supersolvability and solvability of a finite group by its order sequence. We can formulate this problem as follows.\\

\textbf{Question 1.5.} \textit{Let $G$ and $H$ be finite groups. If $G$ is supersolvable/solvable and $\os(G)=\os(H)$, is it true that $H$ is supersolvable/solvable?}\\

The second objective of our paper is to give a negative answer to the 'supersolvable' version of the above question. More comments and detailed proofs are outlined in the following section which deal with these questions posed by Cameron and Dey.
    
\section{Main results}

As a starting point, in the table below, we list all possible pairs $(H, G)$ of finite groups of same order (up to 1000) such that $H$ is non-solvable, $G$ is solvable and $\os(H)$ properly dominates $\os(G)$. The existence of these examples is determined via GAP \cite{13}. They confirm that the answer to Question 1.3 is negative.\\

\noindent\begin{tabular}{ |p{0.8cm}|p{0.2cm}|p{5.85cm}|p{0.35cm}|p{5.8cm}| }
 \hline
 \multicolumn{5}{|c|}{$H\cong SmallGroup(Order, i), G\cong SmallGroup(Order, j)$ and $\os(H)$ properly dominates $\os(G)$}\\
 \hline
$Order$ & $i$ & $\os(H)$  & $j$ & $\os(G)$\\
 \hline
300 & 22 & ((1, 1), (2, 15), (3, 20), (5, 124), (10, 60), (15, 80))
 & 23 & ((1, 1), (2, 25), (3, 50), (4, 150), (5, 24), (6, 50))\\
 \hline
420 & 13 & ((1, 1), (2, 15), (3, 20), (5, 24), (7, 6), (14, 90), (21, 120), (35, 144)) & 16 & ((1, 1), (2, 47), (3, 14), (5, 4), (6, 154), (7, 6), (10, 28), (14, 30), (15, 56), (30, 56),
 (35, 24))\\
\hline
780 & 13 & ((1, 1), (2, 15), (3, 20), (5, 24), (13, 12),  (26, 180), (39, 240), (65, 288))
 & 16 & ((1, 1), (2, 65), (3, 26), (4, 130), (5, 4), (6, 130), (12, 260), (13, 12), (15, 104), (65, 48))\\
\hline
780 & 13 & ((1, 1), (2, 15), (3, 20), (5, 24), (13, 12),  (26, 180), (39, 240), (65, 288)) & 17 & ((1, 1), (2, 65), (3, 26), (4, 130), (5, 4), (6, 130), (12, 260), (13, 12), (15, 104), (65, 48))\\
\hline
780 & 13 & ((1, 1), (2, 15), (3, 20), (5, 24), (13, 12),  (26, 180), (39, 240), (65, 288)) &  20 & ((1, 1), (2, 83), (3, 26), (5, 4), (6, 286), (10, 52), (13, 12), (15, 104), (26, 60), (30, 104), (65, 48))\\
\hline
900 & 88 & ((1, 1), (2, 15), (3, 62), (5, 124), (6, 30), (10, 60), (15, 488), (30, 120)) & 89 & ((1, 1), (2, 25), (3, 152), (4, 150), (5, 24), (6, 200), (12, 300), (15, 48))\\
\cline{4-5}
&  &  & 90 & ((1, 1), (2, 25), (3, 152), (4, 450), (5, 24), (6, 200), (15, 48))\\
\cline{4-5}
 &  &  &  91 & ((1, 1), (2, 25), (3, 152), (4, 50), (5, 24), (6, 200), (12, 400), (15, 48))\\
\cline{4-5}
 &  &  & 92 & ((1, 1), (2, 25), (3, 152), (4, 150), (5, 24), (6, 200), (12, 300), (15, 48))\\
\hline
\end{tabular}

\noindent\begin{tabular}{ |p{0.8cm}|p{0.2cm}|p{5.85cm}|p{0.35cm}|p{5.8cm}| }
 \hline
900 & 88 & ((1, 1), (2, 15), (3, 62), (5, 124), (6, 30), (10, 60), (15, 488), (30, 120)) & 94 & ((1, 1), (2, 63), (3, 152), (5, 24), (6, 180), (10, 312), (15, 48), (30, 120))\\
\cline{4-5}
 &  &  & 95 & ((1, 1), (2, 135), (3, 152), (5, 24), (6, 180), (10, 240), (15, 48), (30, 120))\\
\cline{4-5}
 &  &  & 96 & ((1, 1), (2, 115), (3, 152), (5, 24), (6, 200), (10, 360), (15, 48))\\
\cline{4-5}
 &  &  & 97 & ((1, 1), (2, 103), (3, 152), (5, 24), (6, 500), (10, 72), (15, 48))\\
\cline{4-5}
 &  &  & 97 & ((1, 1), (2, 103), (3, 152), (5, 24), (6, 500), (10, 72), (15, 48))\\
\cline{4-5}
 & &  &  100 & ((1, 1), (2, 91), (3, 152), (5, 24), (6, 152), (10, 384), (15, 48), (30, 48))\\
\cline{4-5}
& &  & 101 & ((1, 1), (2, 51), (3, 152), (5, 24), (6, 552), (10, 24), (15, 48), (30, 48))\\
\cline{4-5}
&  & &  103 & ((1, 1), (2, 151), (3, 152), (5, 24), (6, 452), (10, 24), (15, 48), (30, 48))\\
\cline{4-5}
&  &  & 119 & ((1, 1), (2, 225), (3, 8), (4, 450), (5, 24), (15, 192))\\
\cline{4-5}
 & &  & 120 & ((1, 1), (2, 225), (3, 8), (4, 450), (5, 24), (15, 192))\\
\cline{4-5}
 &  &  &  129 & ((1, 1), (2, 259), (3, 8), (5, 24), (6, 200), (10, 216), (15, 192))\\
\hline
\end{tabular}\\  

The structure descriptions of the groups highlighted on the first row are $H\cong C_5\times A_5$ and $G\cong C_5^2\rtimes Dic_{12}$, respectively. Here, $Dic_{12}$ is the dicyclic group with 12 elements. The structure descriptions of the other 3 non-solvable groups outlined in the table are $C_7\times A_5, C_{13}\times A_5$ and $C_{15}\times A_5$, respectively. Up to isomorphism, these 4 non-solvable groups are the only non-solvable groups of order 300, 420, 780 or 900, respectively. In other words, these 4 numbers are almost solvable numbers (an almost solvable number is a positive integer $n$ such that there is a unique non-solvable group of order $n$). In \cite{11}, the almost abelian and almost nilpotent numbers were fully described, while the almost cyclic numbers were characterized in \cite{3}. The order sequence of $C_{15}\times A_5$ properly dominates the order sequences corresponding to 15 solvable groups. 

We are going to show that we can add infinitely many rows to the previous table by taking $H\cong C_{5p}\times A_5$ and $G\cong C_p\times (C_5^2\rtimes Dic_{12})$, where $p\neq 5$ is an odd prime such that $p\not\equiv 1 \ (mod \ 5)$. After choosing a prime $p$ with the previous properties, an important step is to show that there is only one non-solvable group of order $300p$, it being isomorphic to $H$. This step must be justified in order to avoid the possibility of the existence of a second non-solvable group of order $300p$ whose order sequence is properly dominated by the order sequence of all solvable groups of the same order. 

We denote the Schur multiplier and the outer automorphism group of a finite group $G$ by $M(G)$ and $Out(G)$, respectively. We recall that a simple $K_n$-group is a non-abelian simple group whose order is divisible by exactly $n$ prime numbers. We use the same notations as the ones from \cite{7} for writing various non-abelian simple groups. Our reasoning starts with a lemma which enumerates some results that are going to be further used.\\

\textbf{Lemma 2.1.} \textit{Let $G$ be a finite group and let $S$ be a non-abelian simple group. The following results hold:
\begin{itemize}
\item[i)] If $G$ is non-solvable, then $G$ has a normal series $1\triangleleft H\triangleleft K\triangleleft G$ such that $\frac{K}{H}$ is isomorphic to a direct product of isomorphic non-abelian simple groups and $\big|\frac{G}{K} \big| \big| \big|Out(\frac{K}{H})\big|$ (see Lemma 1 of \cite{14});
\item[ii)] If $G$ is a simple $K_3$-group, then $G$ is isomorphic to $A_5, A_6, L_2(7), L_2(8), L_2(17), L_3(3), U_3(3)$ or $U_4(2)$ (see \cite{9} or Lemma 4 of \cite{8});
\item[iii)] If $G$ is a simple $K_4$-group, then $G$ is isomorphic to $A_7, A_8, A_9, A_{10}, M_{11}, M_{12}, J_2, L_3(4), L_3(5),$ $L_3(7), L_3(8), L_3(17), L_4(3), S_4(4), S_4(5), S_4(7), S_4(9), S_6(2), O^+_8(2), G_2(3), U_3(4), U_3(5), U_3(7),$\\ $U_3(8), U_3(9), U_4(3), U_5(2), Sz(8), Sz(32), ^3D_4(2), ^2F_4(2)'$ or $L_2(q)$, where $q\geq 11$ is a prime power such that $q(q^2-1)=(2, q-1)2^{\alpha_1}3^{\alpha_2}r^{\alpha_3}t^{\alpha_4}$, $\alpha_i>0$ for all $i\in\{1, 2, 3, 4\}$, and $t>r>3$ are primes (see Theorem 1 of \cite{4});
\item[iv)] If $H$ is a cyclic normal subgroup of $G$ such that $\frac{G}{H}\cong S$ and $(|M(S)|, |H|)=1$, then $G\cong H\times S$ (see Theorem 2.1 of \cite{12});
\item[v)] If $H$ is a cyclic normal subgroup of $G$ such that $\frac{G}{H}$ is supersolvable, then $G$ is supersolvable (see Theorem 4.15 of \cite{6}).
\end{itemize}}

\textbf{Theorem 2.2.} \textit{Let $p\neq 5$ be an odd prime such that $p\not\equiv 1 \ (mod \ 5)$ and let $G$ be a non-solvable group of order $300p$. Then $G\cong C_{5p}\times A_5$.}\\

\textbf{Proof.} We already know that $C_{15}\times A_5$ is the only non-solvable group of order 900, so we can assume that $p\geq 7$. By Lemma 2.1, \textit{i)}, we know that $G$ has a normal series $1\triangleleft H\triangleleft K\triangleleft G$ such that $\frac{K}{H}$ is isomorphic to a direct product of isomorphic non-abelian simple groups and $\big|\frac{G}{K} \big| \big| \big|Out(\frac{K}{H})\big|$. Since $|G|=300p$, the section $\frac{K}{H}$ must actually be isomorphic to a non-abelian simple group. This together with the fact that the order of $\frac{K}{H}$ divides $|G|$, leads to
$$\bigg| \frac{K}{H}\bigg|\in \bigg\{60, 12p, 20p, 50p, 60p, 100p, 150p, 300p\bigg\}.$$ Note that none of these potential orders is divisible by 8.

\textit{Case 1.} Assume that $\big|\frac{K}{H}\big|\in\{60, 12p, 20p, 50p, 100p\}$. Then $\frac{K}{H}$ is a simple $K_3$-group. By inspecting the list given by Lemma 2.1, \textit{ii)}, excepting $A_5$, the order of a simple $K_3$-group is divisible by 8. Hence, the only possible choice is $\frac{K}{H}\cong A_5$.

\textit{Case 2.} Assume that $\big|\frac{K}{H}\big|\in\{60p, 150p, 300p\}$. Then $\frac{K}{H}
$ is a simple $K_4$-group, so it must be isomorphic to one of the groups listed in the statement of Lemma 2.1, \textit{iii)}. All possibilities excepting the last one - for an odd prime power $q\geq 11$ - are eliminated since all other simple $K_4$-groups are having  orders divisible by 8. Hence, we get that $\frac{K}{H}\cong L_2(q)$, where $q\geq 11$ is an odd prime power such that $q(q^2-1)=2^{\alpha_1+1}3^{\alpha_2}r^{\alpha_3}t^{\alpha_4}$, $\alpha_i>0$ for all $i\in\{1, 2, 3, 4\}$, and $t>r>3$ are primes. Since $q$ is odd, we have 
$$|L_2(q)|=q\cdot \frac{q^2-1}{2}=2^{\alpha_1}3^{\alpha_2}r^{\alpha_3}t^{\alpha_4}.$$ Also, it is clear that $r=5, t=p, \alpha_2=\alpha_4=1$ and $\alpha_1,\alpha_3\in \{1, 2\}$ as a consequence of the equality $\big|\frac{K}{H}\big|=|L_2(q)|$. As $q\geq 11$ is an odd prime power that divides $|L_2(q)|$, it follows that $q\in \{5^{\alpha_3}, p\}$. If $q=5^{\alpha_3}\geq 11$, then $\alpha_3=2$ and $\big|\frac{K}{H}\big|=|L_2(25)|=7800$. We deduce that $p\in \{26, 52, 130\}$, contradicting that $p$ is a prime. If $q=p$, then $\frac{p^2-1}{2}\in \{ 60, 150, 300\}$. Since $p$ is a prime, we obtain that $p=11$, but this contradicts the initial assumption of working with $p\not\equiv 1 \ (mod \ 5).$

As a consequence of the previous two cases, we infer that $\frac{K}{H}\cong A_5$. Since $|K|$ and $|H|$ are divisors of $G$, we get
$$(|K|, |H|)\in \{ (60, 1), (60p, p), (300, 5), (300p, 5p)\}.$$
Recall that $\big|\frac{G}{K} \big| \big| \big|Out(\frac{K}{H})\big|$. Since $\frac{K}{H}\cong A_5$ and $|Out(A_5)|=2$, the only possible choice is $(|K|, |H|)=(300p, 5p)$. Therefore, $K=G$. This means that $H\triangleleft G$. Moreover, $H\cong C_{5p}$, as $p\not\equiv 1 \ (mod \ 5)$ and $p\geq 7$, $\frac{G}{H}\cong A_5$ and $(|M(A_5)|, |H|)=(2, 15p)=1$. By Lemma 2.1, \textit{iv)}, it follows that $G\cong C_{5p}\times A_5$ and we conclude the proof.    
\hfill\rule{1,5mm}{1,5mm}\\

We are now able to negatively answer to Question 1.3.\\

\textbf{Theorem 2.3.} \textit{There are infinitely many positive integers $n$ such that there exists a unique non-solvable group $H$ of order $n$ whose order sequence properly dominates the order sequence of a solvable group $G$ of the same order.}\\

\textbf{Proof.} By Dirichlet's theorem, we know that there are infinitely many primes $p$ such that $p\not\equiv 1 \ (mod \ 5)$. For each choice of $p\neq 5$, we take $n=300p$. Theorem 2.2 guarantees the existence of a unique non-solvable group $H\cong C_{5p}\times A_5\cong C_p\times (C_5\times A_5)$ of order $n$. Let  $G\cong C_p\times (C_5^2\rtimes Dic_{12})$ be a solvable group of the same order. If $p=3$, we already know that $\os(H)$ properly dominates $\os(G)$ (check the 6th row of the previous table). Assume that $p\geq 7$. Then $(p, 300)=1$. Since $\os(C_5\times A_5)$ properly dominates $\os(C_5^2\rtimes Dic_{12})$, according to Theorem 1.1, \textit{i)}, \textit{iii)}, we conclude that $\os(H)=\os(C_p)\os(C_5\times A_5)$ properly dominates $\os(G)=\os(C_p)\os(C_5^2\rtimes Dic_{12})$.
\hfill\rule{1,5mm}{1,5mm}\\

We mention that an answer to a 'supersolvable' version of Question 1.3 is also negative. Up to order 500, there are several examples of pairs $(H, G)$ of groups of same order such that $H$ is non-supersolvable, $G$ is supersolvable and $\os(H)$ properly dominates $\os(G)$. Without mentioning the order sequences corresponding to these pairs, we resume our findings in the table below.\\

\noindent\begin{tabular}{ |p{1cm}|p{1cm}|p{11.8cm}|  }
 \hline
 \multicolumn{3}{|c|}{$H\cong SmallGroup(Order, i), G\cong SmallGroup(Order, j)$ and $\os(H)$ properly dominates $\os(G)$} \\
 \hline
 \textit{Order} & $i$ & $j$\\
 \hline
 112 & 41 & 36, 42\\
 \hline
 132 & 6 & 5 \\
 \hline
 204 & 8 & 5, 6, 7\\
 \hline
 224 & 173 & 12, 66, 67, 68, 69, 70, 71, 75, 76, 77, 78, 79, 80, 81, 86, 87, 88, 89, 90, 91, 92, 93, 122, 123, 124, 125, 132, 133, 134, 135, 141, 142, 147, 148, 175, 176, 177, 178, 179, 180, 181, 182, 183, 184, 185, 186, 187, 188, 196\\
 \hline
 224 & 195 & 148, 188, 196\\
 \hline
 276 & 6 & 5\\ 
 \hline
 280 & 33 & 32 \\
 \hline
 294 & 7 & 13, 14\\
 \hline
 348 & 8 & 5, 6, 7\\
 \hline
 420 & 13 & 16\\
 \hline
 420 & 32 & 14, 15, 16, 17, 18, 19, 24, 25, 26, 27, 28, 29, 30, 31\\
 \hline
 420 & 33 & 16 \\
 \hline
 492 & 8 & 5, 6, 7\\
 \hline
\end{tabular}\\

As in the solvable case, the 13 non-supersolvable groups outlined in the previous table are the only non-supersolvable groups of those specific orders. In particular, we observe that there are 3 supersolvable groups of order 224, $C_2^4\rtimes D_{14}$ ($j=148$), $C_2^2\times (C_7\rtimes D_8)$ ($j=188$) and $C_2^4\times D_{14}$ ($j=196$), whose order sequences are properly dominated by the order sequences of both non-supersolvable groups of order 224,  $C_4\times F_8$ ($i=173$) and $C_2^2\times F_8$ ($i=195$). Also, there is a supersolvable group of order 420, $D_{10}\times F_7$ ($j=16$) whose order sequence is properly dominated by the order sequences of all 3 non-supersolvable groups of the same order: $C_7\times A_5$ ($i=13$), $C_{35}\times A_4$ ($i=32$), and $C_5\times (C_7\rtimes A_4)$ ($i=33$). Here, $F_7$ and $F_8$ are the Frobenius groups of orders 42 and 56, respectively. We are going to prove that we can add infinitely many examples in the previous table. For an odd prime $p\geq 11$ with $p\not \equiv 1 \ (mod \ 3)$,  the candidates for $H$ and $G$ are $C_p\times A_4$ and $S_3\times D_{2p}$, respectively.\\

\textbf{Theorem 2.4.} \textit{Let $G$ be a non-supersolvable group of order $12p$, where $p\geq 11$ is an odd prime such that $p\not\equiv 1 \ (mod \ 3)$. Then $G\cong C_p\times A_4$.}\\

\textbf{Proof.} If $p=11$, we already know that $G\cong C_{11}\times A_4$ is the only non-supersolvable group of order 132 (SmallGroup(132,6) outlined in the previous table). Assume that $p>11$. It is easy to show that $G$ has a unique Sylow $p$-subgroup $H\cong C_p$. Hence, $H\triangleleft G$. Since $(|H|, \big|\frac{G}{H}\big|)=(p, 12)=1$, by applying the Schur-Zassenahaus theorem, we deduce that $G\cong H\rtimes \frac{G}{H}$. Further, we note that $H$ is a cyclic normal subgroup of a non-supersolvable group $G$. Hence, by Lemma 2.1, \textit{v)}, it follows that $\frac{G}{H}$ is a non-supersolvable group of order 12. Then $\frac{G}{H}\cong A_4$, so $G\cong C_p\rtimes A_4$. Let $f:C_p\longrightarrow Aut(C_p)$ the group homomorphism that induces the previous semidirect product. Let $\sigma\in A_4$ with $o(\sigma)=3$. It is known that $o(f(\sigma))|o(\sigma)$, so $o(f(\sigma))\in\{1, 3\}$. Assume that $o(f(\sigma))=3$. Since $o(f(\sigma))| |Aut(C_p)|$ and $Aut(C_p)\cong C_{p-1}$, it follows that $3|p-1$, a contradiction. Then the trivial permutation and all 8 3-cycles of $A_4$ belong to $\ker f$ whose order must divide 12. Consequently, $\ker f=A_4$, so $f$ is trivial and $G\cong C_p\times A_4$, as desired.     
\hfill\rule{1,5mm}{1,5mm}\\

We have the necessary tools to justify that the answer to the 'supersolvable' version of Question 1.3 is negative.\\

\textbf{Theorem 2.5.} \textit{There are infinitely many positive integers $n$ such that there exists a unique non-supersolvable group $H$ of order $n$ whose order sequence properly dominates the order sequence of a supersolvable group $G$ of the same order.}\\

\textbf{Proof.} By Dirichlet's theorem, we know that there are infinitely many odd primes $p\geq 11$ such that $p\equiv 2 \ (mod \ 3)$. For each such $p$, we take $n=12p$. Theorem 2.5 assures the existence of a unique non-supersolvable group $H\cong C_{p}\times A_4$ of order $n$. Let $G\cong S_3\times D_{2p}$ be a supersolvable group of the same order. We are going to effectively write the order sequences of $H$ and $G$ and compare them. We recall that $os(A_4)=(1, 2, 2, 2, 3, 3, 3, 3, 3, 3, 3, 3)$. Hence,
$$\os(H)=\os(C_p\times A_4)=((1,1), (2,3), (3,8), (p, p-1), (2p, 3p-3), (3p, 8p-8)).$$
The order sequences of $S_3$, $D_{2p}$ are $\os(S_3)=(1, 2, 2, 2, 3, 3)$ and $\os(D_{2p})=((1,1), (2, p), (p, p-1))$, respectively. Then,
$$\os(G)=\os(S_3\times D_{2p})=((1, 1), (2, 4p+3), (3, 2), (6, 2p), (p, p-1), (2p, 3p-3), (3p, 2p-2)).$$ 
By comparing $\os(H)$ with $\os(G)$, we observe that the first sequence properly dominates the second one and this concludes the proof.
\hfill\rule{1,5mm}{1,5mm}\\

As a consequence of Theorems 2.2 and 2.4, we find infinitely many almost solvable/supersolvable numbers.\\

\textbf{Corollary 2.6.} \textit{\begin{itemize}
\item[i)] The numbers $300p$, where $p\neq 5$ is an odd prime satisfying $p\not\equiv 1 \ (mod \ 5)$, are almost solvable numbers.
\item[ii)] The numbers $12p$, where $p\geq 11$ is an odd prime satisfying $p\not \equiv 1 \ (mod \ 3)$, are almost supersolvable numbers.
\end{itemize}}

Before moving our attention to the 'supersolvable' version of Question 1.5, we pose some open problems that could further develop our findings. It can be checked that $\os(A_5)$, $\os(L_2(7))$, $\os(A_6)$, $\os(L_2(8))$, $\os(L_2(11))$ and $\os(L_2(13))$ are minimal elements in their corresponding posets. In \cite{10}, it was proved that $\psi(C_3^2\times Sz(8))<\psi(L_2(64))$. This showed that a conjecture (see \cite{1}) stating that if a finite simple group $S$ and a non-simple group $G$ are of same order, then $\psi(S)<\psi(G)$, is false. However, the order sequences of these two groups are incomparable with respect to the domination relation since
\begin{align*}
\os(L_2(64))=&((1, 1), (2, 4095), (3, 4160), (5, 8064), (7, 12480), (9, 12480), (13, 24192), (21, 24960),\\&
  (63, 74880), (65, 96768))
\end{align*}
and
\begin{align*}
\os(C_3^2\times Sz(8))=&((1, 1), (2, 455), (3, 8), (4, 3640),  (5, 5824), (6, 3640), (7, 12480), (12, 29120),\\& (13, 6720),
  (15, 46592), (21, 99840), (39, 53760)).
\end{align*}

\textbf{Question 2.7.} \textit{Let $n$ be a positive integer such that there exists a simple group $G$ of order $n$. Is is true that $\os(G)$ is a minimal element in the poset formed of the order sequences of all groups of order $n$ with respect to the domination relation?}\\

We determined a class of almost solvable numbers, their form being $300p$ where $p\neq 5$ is an odd prime with $p\not \equiv 1 \ (mod \ 5)$. Also, we indicated a class of almost supersolvable numbers of the form $12p$ with $p\geq 11$ being an odd prime such that $p\not\equiv 1 \ (mod \ 3)$. As we previously highlighted, 420 and 780 are other almost solvable numbers, while 112, 280 and 294 are almost supersolvable numbers.\\

\textbf{Question 2.8.} \textit{Describe the almost supersolvable/solvable numbers.}\\

Finally, we are going to negatively answer to the 'supersolvable' version of Question 1.5. As a consequence, we can state that the supersolvability of a finite group cannot be determined by its order sequence. Up to order 250, there are several pairs $(H, G)$ of groups of the same order order such that $H$ is non-supersolvable, $G$ is supersolvable and $\os(G)=\os(H)$. We list them below.\\   

\noindent\begin{tabular}{ |p{1cm}|p{1.3cm}|p{1.6cm}|p{9.6cm}|  }
 \hline
 \multicolumn{4}{|c|}{$H\cong SmallGroup(Order, i), G\cong SmallGroup(Order, j)$ and $\os(G)=\os(H)$} \\
 \hline
 \textit{Order} & $i$ & $j$ & $\os(H)$\\
 \hline
 72 & 40 & 35 & ((1, 1), (2, 21), (3, 8), (4, 18), (6, 24))\\
 \hline
 144 & 119 & 99 & ((1, 1), (2, 1), (3, 8), (4, 42), (6, 8), ( 8, 36), (12, 48))\\
 \hline
144 & 116 & 174 & ((1, 1), (2, 19), (3, 8), (4, 60), (6, 8), (12, 48))\\
\hline
144 & 186 & 177 & ((1, 1), (2, 43), (3, 8), (4, 36), (6, 56))\\
\hline
 216 & 100, 168 & 34, 36, 131 & ((1, 1), (2, 19), (3, 26), (4, 36), (6, 62), (12, 72))\\
 \hline
 216 & 157 & 60, 72, 144 & ((1, 1), (2, 21), (3, 26), (4, 18), (6, 114), (12, 36))\\
 \hline
\end{tabular}\\

We observe that there are 3 supersolvable groups and 2 non-supersolvable groups of order 216 having equal order sequences. The structure descriptions of the groups outlined on the first row are $H\cong S_3^2\rtimes C_2$ and $G\cong C_3^2\rtimes D_8$.  By starting with this first pair $(H, G)$, we can build infinitely many other examples.\\

\textbf{Theorem 2.9.} \textit{Let $p\geq 5$ be a prime number. Then there exists a pair $(H, G)$ of groups of order $72p$ such that $H$ is non-supersolvable, $G$ is supersolvable and $\os(G)=\os(H)$.}\\

\textbf{Proof.} Let  $p\geq 5$ be a prime. Let $H\cong C_p\times (S_3^2\rtimes C_2)$ be a non-supersolvable group of order $72p$ and let $G\cong C_p\times (C_3^2\rtimes D_8)$ be a supersolvable group of the same order. By Theorem 1.1, \textit{i)}, and by using the equality $\os(S_3^2\rtimes C_2)=\os(C_3^2\rtimes D_8)$, we have
$$\os(G)=\os(C_p)\os(C_3^2\rtimes D_8)=\os(C_p)\os(S_3^2\rtimes C_2)=\os(H),$$
as desired.
\hfill\rule{1,5mm}{1,5mm}\\ 

In the end, we mention that the answer to the 'solvable' version of Question 1.5 is affirmative up to order 2015.

\vspace*{3ex}
\small
\hfill
\begin{minipage}[t]{7cm}
Mihai-Silviu Lazorec \\
Faculty of  Mathematics \\
"Al.I. Cuza" University \\
Ia\c si, Romania \\
e-mail: {\tt silviu.lazorec@uaic.ro}
\end{minipage}
\end{document}